\documentclass[12pt]{article}
\usepackage{amsmath,amssymb,amsthm}

\newtheorem{theorem}{Theorem}
\newtheorem{proposition}{Proposition}
\newtheorem{lemma}[theorem]{Lemma}

\begin{document}

\title{Multivector Functions of a Multivector Variable\thanks{%
published: \emph{Advances in Applied Clifford Algebras} \textbf{11}(S3),
79-91 (2001).}}
\author{A. M. Moya$^{1}\thanks{%
e-mail: moya@ime.unicamp.br}$, V. V. Fern\'{a}ndez$^{1}\thanks{%
e-mail: vvf@ime.unicamp.br}$ and W. A. Rodrigues Jr.$^{1,2}\thanks{%
e-mail: walrod@ime.unicamp.br or walrod@mpc.com.br}$ \\
$\hspace{-0.5cm}^{1}$ Institute of Mathematics, Statistics and Scientific
Computation\\
IMECC-UNICAMP CP 6065\\
13083-970 Campinas-SP, Brazil\\
$^{2}$ Department of Mathematical Sciences, University of Liverpool\\
Liverpool, L69 3BX, UK}
\date{11/26/2001}
\maketitle

\begin{abstract}
In this paper we develop with considerable details a theory of multivector
functions of a $p$-vector variable. The concepts of limit, continuity and
differentiability are rigorously studied. Several important types of
derivatives for these multivector functions are introduced, as e.g., the $A$%
-directional derivative (where $A$ is a $p$-vector) and the generalized
concepts of curl, divergence and gradient. The derivation rules for
different types of products of multivector functions and for compositon of
multivector functions are proved.
\end{abstract}

\tableofcontents

\section{Introduction}

This is the paper VI in the present series. Here, we develop a theory of
multivector functions of a $p$-vector variable. For these objects we
investigate with details the concepts of limit and continuity, and formulate
rigorously the notion of derivation. As we will see, the concept of extensor
introduced in \cite{1} (paper II on this series) plays a crucial role in our
theory of differentiability. We introduce important derivative-like
operators for these multivector functions, as e.g., the $A$-directional
derivative and the generalized concepts of curl, divergence and gradient.
The derivation rules for all suitable products of multivector functions of a
$p$-vector variable and for composition of multivector functions are
presented and proved.

\section{Multivector Functions of a $p$-Vector Variable}

Let $\mathbf{\Omega}^{p}V$ be a subset of $\bigwedge^{p}V.$ Any mapping $F:%
\mathbf{\Omega}^{p}V\rightarrow\bigwedge V$ will be called a \emph{%
multivector function of a }$p$\emph{-vector variable over }$V.$ In
particular, $F:\mathbf{\Omega}^{p}V\rightarrow\bigwedge^{q}V$ is said to be
a $q$\emph{-vector function of a }$p$\emph{-vector variable, }or a $(p,q)$%
\emph{-function over }$V,$ for short. For the special cases $q=0,$ $q=1,$ $%
q=2,\ldots$ etc. we will employ the names of \emph{scalar, vector, bivector,}%
$\ldots$\emph{\ etc. function of a }$p$\emph{-vector variable,} respectively.

\subsection{Limit Notion}

We begin by introducing the concept of $\delta$-neighborhood for a
multivector $A.$

Take any $\delta>0.$ The set\footnote{%
We recalls that the two double bars $\left\| \left. {}\right. \right\| $
denotes the norm of multivectors, as defined in [2], i.e., for all $%
X\in\Lambda V:\left\| X\right\| =\sqrt{X\cdot X},$ where $(\cdot)$ is any
fixed euclidean scalar product.} $N_{A}(\delta)=\{X\in\Lambda V$ / $\left\|
X-A\right\| <\delta\}$ will be called a $\delta$\emph{-neighborhood of }$A. $

The set $N_{A}(\delta)-\{A\}=\{X\in\bigwedge V$ / $0<\left\| X-A\right\|
<\delta\}$ will be said to be a \emph{reduced }$\delta$\emph{-neighborhood
of }$A.$

We introduce now the concepts of cluster and interior points of $\mathbf{%
\Omega}V\subseteq\bigwedge V$.

A multivector $X_{0}\in\bigwedge V$ is said to be a \emph{cluster point of }$%
\mathbf{\Omega}V$ if and only if for every $N_{X_{0}}(\delta):(N_{X_{0}}(%
\delta)-\{X_{0}\})\cap\mathbf{\Omega}V\neq\emptyset,$ i.e., all reduced $%
\delta$-neighborhood of $X_{0}$ contains at least one multivector of $%
\mathbf{\Omega}V.$

A multivector $X_{0}\in\bigwedge V$ is said to be an \emph{interior point of
}$\mathbf{\Omega}V$ if and only if there exists $N_{X_{0}}(\delta)$ such
that $N_{X_{0}}(\delta)\subseteq\mathbf{\Omega}V,$ i.e., any multivector of
some $\delta$-neighborhood of $X_{0}$ belongs also to $\mathbf{\Omega}V.$

It should be noted that any interior point of $\mathbf{\Omega}V$ is also a
cluster point of $\mathbf{\Omega}V$.

If the set of interior points of $\mathbf{\Omega}V$ coincides with $\mathbf{%
\Omega}V,$ then it is said to be an \emph{open subset of }$\bigwedge V.$

Take $\mathbf{\Omega}^{p}V\subseteq\bigwedge^{p}V$ and let $F:\mathbf{\Omega
}^{p}V\rightarrow\bigwedge V$ be a multivector function of a $p$-vector
variable and take $X_{0}\in\bigwedge^{p}V$ to be a cluster point of $\mathbf{%
\Omega}^{p}V.$

A multivector $M$ is said to be the \emph{limit of }$F(X)$\emph{\ for }$X$%
\emph{\ approaching to }$X_{0}$ if and only if for every real $\varepsilon>0
$ there exists some real $\delta>0$ such that if $X\in \mathbf{\Omega}^{p}V $
and $0<\left\| X-X_{0}\right\| <\delta,$ then $\left\| F(X)-M\right\|
<\varepsilon.$ It is denoted by $\underset{X\rightarrow X_{0}}{\lim}F(X)=M.$

In dealing with a scalar function of a $p$-vector variable, say $\Phi,$ the
definition of $\underset{X\rightarrow X_{0}}{\lim}\Phi(X)=\mu$ is reduced
to: for every $\varepsilon>0$ there exists some $\delta>0$ such that $\left|
\Phi(X)-\mu\right| <\varepsilon,$ whenever $X\in\Omega^{p}V $ and $0<\left\|
X-X_{0}\right\| <\delta.$

\begin{proposition}
Let $F:\mathbf{\Omega}^{p}V\rightarrow\bigwedge V$ and $G:\mathbf{\Omega}%
^{p}V\rightarrow\bigwedge V$ be two multivector functions of a $p$-vector
variable. If there exist $\underset{X\rightarrow X_{0}}{\lim}F(X)$ and $%
\underset{X\rightarrow X_{0}}{\lim}G(X),$ then there exists $\underset{%
X\rightarrow X_{0}}{\lim}(F+G)(X)$ and
\begin{equation}
\underset{X\rightarrow X_{0}}{\lim}(F+G)(X)=\underset{X\rightarrow X_{0}}{%
\lim}F(X)+\underset{X\rightarrow X_{0}}{\lim}G(X).   \label{6.1a}
\end{equation}
\end{proposition}

\begin{proof}
Let $\underset{X\rightarrow X_{0}}{\lim}F(X)=M_{1}$ and $\underset{%
X\rightarrow X_{0}}{\lim}G(X)=M_{2}.$ Then, we must prove that $\underset{%
X\rightarrow X_{0}}{\lim}(F+G)(X)=M_{1}+M_{2}.$

Taken a real $\varepsilon>0,$ since $\underset{X\rightarrow X_{0}}{\lim }%
F(X)=M_{1}$ and $\underset{X\rightarrow X_{0}}{\lim}G(X)=M_{2},$ there must
be two real numbers $\delta_{1}>0$ and $\delta_{2}>0$ such that
\begin{align*}
\left\| F(X)-M_{1}\right\| & <\frac{\varepsilon}{2},\text{ for }X\in
\Omega^{p}V\text{ and }0<\left\| X-X_{0}\right\| <\delta_{1}, \\
\left\| G(X)-M_{2}\right\| & <\frac{\varepsilon}{2},\text{ for }X\in
\Omega^{p}V\text{ and }0<\left\| X-X_{0}\right\| <\delta_{2}.
\end{align*}

Thus, there is a real $\delta=\min\{\delta_{1},\delta_{2}\}$ such that
\begin{equation*}
\left\| F(X)-M_{1}\right\| <\frac{\varepsilon}{2}\text{ and }\left\|
G(X)-M_{2}\right\| <\frac{\varepsilon}{2},
\end{equation*}
for $X\in\mathbf{\Omega}^{p}V$ and $0<\left\| X-X_{0}\right\| <\delta.$
Hence, by using the triangular inequality for the norm of multivectors, it
follows that
\begin{align*}
\left\| (F+G)(X)-(M_{1}+M_{2})\right\| & =\left\|
F(X)-M_{1}+G(X)-M_{2}\right\| \\
& \leq\left\| F(X)-M_{1}\right\| +\left\| G(X)-M_{2}\right\| \\
& <\frac{\varepsilon}{2}+\frac{\varepsilon}{2}=\varepsilon,
\end{align*}
for $X\in\mathbf{\Omega}^{p}V$ and $0<\left\| X-X_{0}\right\| <\delta.$

Therefore, for any $\varepsilon>0$ there is a $\delta>0$ such that if $X\in%
\mathbf{\Omega}^{p}V$ and $0<\left\| X-X_{0}\right\| <\delta,$ then $\left\|
(F+G)(X)-(M_{1}+M_{2})\right\| <\varepsilon.$
\end{proof}

\begin{proposition}
Let $\Phi:\mathbf{\Omega}^{p}V\rightarrow\mathbb{R}$ and $F:\mathbf{\Omega }%
^{p}V\rightarrow\bigwedge V$ be a scalar function and a multivector function
of a $p$-vector variable. If there exist $\underset{X\rightarrow X_{0}}{\lim
}\Phi(X)$ and $\underset{X\rightarrow X_{0}}{\lim}F(X),$ then there exists $%
\underset{X\rightarrow X_{0}}{\lim}(\Phi F)(X)$ and
\begin{equation}
\underset{X\rightarrow X_{0}}{\lim}(\Phi F)(X)=\underset{X\rightarrow X_{0}}{%
\lim}\Phi(X)\underset{X\rightarrow X_{0}}{\lim}F(X).   \label{6.1b}
\end{equation}
\end{proposition}

\begin{proof}
Let $\underset{X\rightarrow X_{0}}{\lim}\Phi(X)=\Phi_{0}$ and $\underset{%
X\rightarrow X_{0}}{\lim}F(X)=F_{0}.$ Then, we must prove that $\underset{%
X\rightarrow X_{0}}{\lim}(\Phi F)(X)=\Phi_{0}F_{0}.$

First, since $\underset{X\rightarrow X_{0}}{\lim}\Phi(X)=\Phi_{0}$ it can be
found a $\delta_{1}>0$ such that
\begin{equation*}
\left| \Phi(X)-\Phi_{0}\right| <1,\text{ for }X\in\mathbf{\Omega}^{p}V\text{
and }0<\left\| X-X_{0}\right\| <\delta_{1},
\end{equation*}
i.e.,
\begin{equation*}
\left| \Phi(X)\right| <1+\left| \Phi_{0}\right| ,\text{ for }X\in\mathbf{%
\Omega}^{p}V\text{ and }0<\left\| X-X_{0}\right\| <\delta_{1}.
\end{equation*}
Where the triangular inequality for real numbers $\left| \alpha\right|
-\left| \beta\right| \leq\left| \alpha-\beta\right| $ was used.

Now, taken a $\varepsilon>0,$ since $\underset{X\rightarrow X_{0}}{\lim}%
\Phi(X)=\Phi_{0}$ and $\underset{X\rightarrow X_{0}}{\lim}F(X)=F_{0},$ they
can be found a $\delta_{2}>0$ and a $\delta_{3}>0$ such that
\begin{align*}
\left| \Phi(X)-\Phi_{0}\right| & <\frac{\varepsilon}{2(1+\left\|
F_{0}\right\| )},\text{ for }X\in\mathbf{\Omega}^{p}V\text{ and }0<\left\|
X-X_{0}\right\| <\delta_{2}, \\
\left\| F(X)-F_{0}\right\| & <\frac{\varepsilon}{2(1+\left| \Phi _{0}\right|
)},\text{ for }X\in\mathbf{\Omega}^{p}V\text{ and }0<\left\| X-X_{0}\right\|
<\delta_{3}.
\end{align*}

Thus, given a real $\varepsilon>0$ there is a real $\delta=\min\{\delta
_{1},\delta_{2},\delta_{3}\}$ such that
\begin{align*}
\left| \Phi(X)\right| & <1+\left| \Phi_{0}\right| , \\
\left| \Phi(X)-\Phi_{0}\right| & <\frac{\varepsilon}{2(1+\left\|
F_{0}\right\| )}, \\
\left\| F(X)-F_{0}\right\| & <\frac{\varepsilon}{2(1+\left| \Phi _{0}\right|
)},
\end{align*}
whenever $X\in\mathbf{\Omega}^{p}V$ and $0<\left\| X-X_{0}\right\| <\delta.$
Hence, using some properties of the norm of multivectors, it follows that
\begin{align*}
\left\| (\Phi F)(X)-\Phi_{0}F_{0}\right\| & =\left\|
\Phi(X)(F(X)-F_{0})+(\Phi(X)-\Phi_{0})F_{0}\right\| \\
& \leq\left| \Phi(X)\right| \left\| F(X)-F_{0}\right\| +\left|
\Phi(X)-\Phi_{0}\right| \left\| F_{0}\right\| \\
& <\left| \Phi(X)\right| \left\| F(X)-F_{0}\right\| +\left| \Phi
(X)-\Phi_{0}\right| (1+\left\| F_{0}\right\| ) \\
& <(1+\left| \Phi_{0}\right| )\frac{\varepsilon}{2(1+\left| \Phi _{0}\right|
)}+\frac{\varepsilon}{2(1+\left\| F_{0}\right\| )}(1+\left\| F_{0}\right\|
)=\varepsilon,
\end{align*}
whenever $X\in\mathbf{\Omega}^{p}V$ and $0<\left\| X-X_{0}\right\| <\delta.$

Therefore, for any $\varepsilon>0$ there is a $\delta>0$ such that if $X\in%
\mathbf{\Omega}^{p}V$ and $0<\left\| X-X_{0}\right\| <\delta,$ then $\left\|
(\Phi F)(X)-\Phi_{0}F_{0}\right\| <\varepsilon.$
\end{proof}

\begin{lemma}
There exists $\underset{X\rightarrow X_{0}}{\lim}F(X)$ if and only if there
exist either $\underset{X\rightarrow X_{0}}{\lim}F^{J}(X)$ or $\underset{%
X\rightarrow X_{0}}{\lim}F_{J}(X).$ It holds
\begin{equation}
\underset{X\rightarrow X_{0}}{\lim}F(X)=\underset{J}{\sum}\frac{1}{\nu (J)!}%
\underset{X\rightarrow X_{0}}{\lim}F^{J}(X)e_{J}=\underset{J}{\sum}\frac{1}{%
\nu(J)!}\underset{X\rightarrow X_{0}}{\lim}F_{J}(X)e^{J}.   \label{6.1c}
\end{equation}
\end{lemma}

\begin{proof}
It is an immediate consequence of eqs.(\ref{6.1a}) and (\ref{6.1b}).
\end{proof}

\begin{proposition}
Let $F:\mathbf{\Omega}^{p}V\rightarrow\bigwedge V$ and $G:\mathbf{\Omega}%
^{p}V\rightarrow\bigwedge V$ be two multivector functions of a $p$-vector
variable. We can define the products $F\ast G:\mathbf{\Omega}%
^{p}V\rightarrow\bigwedge V$ such that $(F\ast G)(X)=F(X)\ast G(X)$ where $%
\ast$ holds for either $(\wedge),$ $(\cdot),$ $(\lrcorner\llcorner)$ or $($%
\emph{Clifford product}$).$ If there exist $\underset{X\rightarrow X_{0}}{%
\lim}F(X)$ and $\underset{X\rightarrow X_{0}}{\lim}G(X),$ then there exists $%
\underset{X\rightarrow X_{0}}{\lim}(F\ast G)(X)$ and
\begin{equation}
\underset{X\rightarrow X_{0}}{\lim}(F\ast G)(X)=\underset{X\rightarrow X_{0}}%
{\lim}F(X)\ast\underset{X\rightarrow X_{0}}{\lim}G(X).   \label{6.1d}
\end{equation}
\end{proposition}

\begin{proof}
It is an immediate consequence of eq.(\ref{6.1c}).
\end{proof}

\subsection{Continuity Notion}

Take $\mathbf{\Omega}^{p}V\subseteq\bigwedge^{p}V.$ A multivector function
of a $p$-vector variable $F:\mathbf{\Omega}^{p}V\rightarrow\bigwedge V$ is
said to be \emph{continuous at }$X_{0}\in\mathbf{\Omega}^{p}V$ if and only
if there exists\footnote{%
Observe that $X_{0}$ has to be cluster point of $\mathbf{\Omega}^{p}V$.} $%
\underset{X\rightarrow X_{0}}{\lim}F(X)$ and
\begin{equation}
\underset{X\rightarrow X_{0}}{\lim}F(X)=F(X_{0}).   \label{6.1e}
\end{equation}

\begin{lemma}
The multivector function $X\mapsto F(X)$ is continuous at $X_{0}$ if and
only if any component scalar function, either $X\mapsto F^{J}(X)$ or $%
X\mapsto F_{J}(X)$ is continuous at $X_{0}$.
\end{lemma}

\begin{proposition}
Let $F:\mathbf{\Omega}^{p}V\rightarrow\bigwedge V$ and $G:\mathbf{\Omega}%
^{p}V\rightarrow\bigwedge V$ be two continuous functions at $X_{0}\in\mathbf{%
\Omega}^{p}V.$ Then, the addition $F+G:\mathbf{\Omega}^{p}V\rightarrow%
\bigwedge V$ such that $(F+G)(X)=F(X)+G(X)$ and the products $F\ast G:%
\mathbf{\Omega}^{p}V\rightarrow\mathbf{\Omega}V$ such that $(F\ast
G)(X)=F(X)\ast G(X),$ where $\ast$ means either $(\wedge),$ $(\cdot),$ $%
(\lrcorner\llcorner)$ or $($\emph{Clifford product}$),$ are also continuous
functions at $X_{0}.$
\end{proposition}

\begin{proof}
It is an immediate consequence of eqs.(\ref{6.1a}) and (\ref{6.1d}).
\end{proof}

\begin{proposition}
Let $G:\mathbf{\Omega}^{p}V\rightarrow\bigwedge^{q}V$ and $F:\bigwedge
^{q}V\rightarrow\bigwedge^{r}V$ be two continuous functions, the first one
at $X_{0}\in\mathbf{\Omega}^{p}V$ and the second one at $G(X_{0})\in%
\bigwedge ^{q}V.$ Then, the composition $F\circ G:\mathbf{\Omega}%
^{p}V\rightarrow \bigwedge^{r}V$ such that $F\circ G(X)=F(G(X))$ is a
continuous function at $X_{0}.$
\end{proposition}

\subsection{Differentiability Notion}

Let $\mathbf{\Omega}^{p}V$ be a subset of $\bigwedge^{p}V.$ A $(p,q)$%
-function over $V,$ say $F,$ is said to be \emph{differentiable at }$X_{0}\in%
\mathbf{\Omega}^{p}V$ if and only if there exists a $(p,q)$-extensor over $V,
$ say $f_{X_{0}},$ such that
\begin{equation}
\underset{X\rightarrow X_{0}}{\lim}\dfrac{F(X)-F(X_{0})-f_{X_{0}}(X-X_{0})}{%
\left\| X-X_{0}\right\| }=0,   \label{6.2a}
\end{equation}
i.e.,
\begin{equation}
\underset{H\rightarrow0}{\lim}\dfrac{F(X_{0}+H)-F(X_{0})-f_{X_{0}}(H)}{%
\left\| H\right\| }=0.   \label{6.2b}
\end{equation}

It is remarkable that if there is such a $(p,q)$-extensor $f_{X_{0}},$ then
it must be \emph{unique}.

Indeed, assume that there is another $(p,q)$-extensor $\overset{\smallfrown }%
{f}_{X_{0}}$ which satisfies
\begin{equation*}
\underset{H\rightarrow0}{\lim}\dfrac{F(X_{0}+H)-F(X_{0})-\overset{%
\smallfrown }{f}_{X_{0}}(H)}{\left\| H\right\| }=0,
\end{equation*}
or equivalently,
\begin{equation*}
\underset{H\rightarrow0}{\lim}\dfrac{\left\| F(X_{0}+H)-F(X_{0})-\overset{%
\smallfrown}{f}_{X_{0}}(H)\right\| }{\left\| H\right\| }=0.
\end{equation*}

By using the triangular inequality which is valid for the norm of
multivectors [1], it can be easily establish the following inequality
\begin{align*}
0 & \leq\frac{\left\| f_{X_{0}}(H)-\overset{\smallfrown}{f}%
_{X_{0}}(H)\right\| }{\left\| H\right\| }\leq\frac{\left\|
F(X_{0}+H)-F(X_{0})-f_{X_{0}}(H)\right\| }{\left\| H\right\| } \\
& +\dfrac{\left\| F(X_{0}+H)-F(X_{0})-\overset{\smallfrown}{f}%
_{X_{0}}(H)\right\| }{\left\| H\right\| },
\end{align*}
which holds for all $H\neq0$ (i.e., $X\neq X_{0}$).

Now, taking the limits for $H\rightarrow0$ (i.e., $X\rightarrow X_{0}$) of
these scalar-valued functions of a $p$-vector variable, we get
\begin{equation*}
\underset{H\rightarrow0}{\lim}\frac{\left\| f_{X_{0}}(H)-\overset{\smallfrown%
}{f}_{X_{0}}(H)\right\| }{\left\| H\right\| }=0.
\end{equation*}

This implies\footnote{%
In order to see that, we can use a lemma: if $\underset{H\rightarrow0}{\lim}%
\Phi(H)=0,$ then $\underset{\lambda \rightarrow0}{\lim}\Phi(\lambda A)=0,$
for all $A\neq0.$} that for every $A\neq0$
\begin{equation*}
\underset{\lambda\rightarrow0}{\lim}\frac{\left\| f_{X_{0}}(\lambda A)-%
\overset{\smallfrown}{f}_{X_{0}}(\lambda A)\right\| }{\left\| \lambda
A\right\| }=0.
\end{equation*}
Then, it follows that for every $A\neq0$
\begin{equation*}
\frac{\left\| f_{X_{0}}(A)-\overset{\smallfrown}{f}_{X_{0}}(A)\right\| }{%
\left\| A\right\| }=0,
\end{equation*}
i.e., $f_{X_{0}}(A)=\overset{\smallfrown}{f}_{X_{0}}(A).$ Now, for $A=0$
this equality trivially holds. Therefore, we have proved that $f_{X_{0}}=%
\overset{\smallfrown}{f}_{X_{0}}$.

The $(p,q)$-extensor $f_{X_{0}}$ will be called the \emph{differential of
the }$(p,q)$\emph{-function }$F$\emph{\ at }$X_{0}.$

So that, the differentiability of $F$ at $X_{0}\in\mathbf{\Omega}^{p}V$
implies the existence of differential of $F$ at $X_{0}\in\mathbf{\Omega}%
^{p}V.$

\begin{lemma}
Associated to any $(p,q)$-function $F,$ differentiable at $X_{0},$ there
exists a $(p,q)$-function $\varphi_{X_{0}},$ continuous at $X_{0,}$ such
that
\begin{equation}
\varphi_{X_{0}}(X_{0})=0   \label{6.2c}
\end{equation}
and for every $X\in\mathbf{\Omega}^{p}V$ it holds
\begin{equation}
F(X)=F(X_{0})+f_{X_{0}}(X-X_{0})+\left\| X-X_{0}\right\| \varphi_{X_{0}}(X).
\label{6.2d}
\end{equation}
\end{lemma}

\begin{proof}
Since the $(p,q)$-function $F$ is differentiable at $X_{0},$ we can define
an auxiliary $(p,q)$-function $\varphi_{X_{0}}$ by
\begin{equation*}
\varphi_{X_{0}}(X)=\left\{
\begin{array}{cc}
0 & \text{for }X=X_{0} \\
\dfrac{F(X)-F(X_{0})-f_{X_{0}}(X-X_{0})}{\left\| X-X_{0}\right\| } & \text{%
for }X\neq X_{0}%
\end{array}
\right. .
\end{equation*}

It satisfies $\varphi_{X_{0}}(X_{0})=0$ and, by taking limit of $\varphi
_{X_{0}}(X)$ for $X\rightarrow X_{0},$ we have
\begin{equation*}
\underset{X\rightarrow X_{0}}{\lim}\varphi_{X_{0}}(X)=\underset{X\rightarrow
X_{0}}{\lim}\frac{F(X)-F(X_{0})-f_{X_{0}}(X-X_{0})}{\left\| X-X_{0}\right\| }%
=0.
\end{equation*}
It follows that $\varphi_{X_{0}}$ is continuous at $X_{0}$ and $\varphi
_{X_{0}}(X_{0})=0.$

Recall now that for $X\neq X_{0}$ it follows the multivector identity
\begin{equation*}
F(X)=F(X_{0})+f_{X_{0}}(X-X_{0})+\left\| X-X_{0}\right\| \varphi_{X_{0}}(X),
\end{equation*}
which for $X=X_{0}$ it is trivially true.
\end{proof}

As happens in the $\mathbb{R}^{n}$ calculus, differentiability implies
continuity. Indeed, by taking limits for $X\rightarrow X_{0}$ on both sides
of eq.(\ref{6.2d}), we get $\underset{X\rightarrow X_{0}}{\lim}F(X)=F(X_{0}).
$

\subsubsection{Directional Derivative}

Since $\mathbf{\Omega}^{p}V$ is an open subset of $\bigwedge^{p}V,$ any $p$%
-vector $X_{0}$ belonging to $\mathbf{\Omega}^{p}V$ is an interior point of $%
\mathbf{\Omega}^{p}V,$ i.e., there is some $\varepsilon$-neighborhood of $%
X_{0},$ say $N_{X_{0}}^{p}(\varepsilon),$ such that $N_{X_{0}}^{p}(%
\varepsilon)\subseteq\mathbf{\Omega}^{p}V.$

Now, take a non-zero $p$-vector $A$ and choose a real number $\lambda$ such
that $0<\left| \lambda\right| <\dfrac{\varepsilon}{\left\| A\right\| }.$
Then, from the obvious inequality $\left\| (X_{0}+\lambda A)-X_{0}\right\|
=\left| \lambda\right| \left\| A\right\| <\varepsilon$ it follows that $%
(X_{0}+\lambda A)\in N_{X_{0}}^{p}(\varepsilon).$ Thus, $(X_{0}+\lambda A)\in%
\mathbf{\Omega}^{p}V.$

There exists $\underset{\lambda\rightarrow0}{\lim}\dfrac{F(X_{0}+\lambda
A)-F(X_{0})}{\lambda}$ and it equals $f_{X_{0}}(A).$

Indeed, by using eq.(\ref{6.2d}) we have
\begin{align*}
\frac{F(X_{0}+\lambda A)-F(X_{0})}{\lambda} & =\frac{f_{X_{0}}(\lambda
A)+\left\| \lambda A\right\| \varphi_{X_{0}}(X_{0}+\lambda A)}{\lambda} \\
& =f_{X_{0}}(A)\pm\left\| A\right\| \varphi_{X_{0}}(X_{0}+\lambda A).
\end{align*}

Now, by taking limits for $\lambda\rightarrow0$ on these $q$-vector
functions of a real variable, using the continuity of $\varphi_{X_{0}}$ at $%
X_{0}$ and eq.(\ref{6.2c}), the required result follows.

The $A$\emph{-directional derivative of }$F$\emph{\ at }$X_{0},$
conveniently denoted by $F_{A}^{\prime}(X_{0}),$ is defined to be
\begin{equation}
F_{A}^{\prime}(X_{0})=\underset{\lambda\rightarrow0}{\lim}\frac{%
F(X_{0}+\lambda A)-F(X_{0})}{\lambda},   \label{6.3a}
\end{equation}

i.e.,
\begin{equation}
F_{A}^{\prime}(X_{0})=\left. \dfrac{d}{d\lambda}F(X_{0}+\lambda A)\right|
_{\lambda=0}.   \label{6.3b}
\end{equation}

The above observation yields a noticeable multivector identity,
\begin{equation}
F_{A}^{\prime}(X_{0})=f_{X_{0}}(A)   \label{6.3c}
\end{equation}
which relates the $A$-directional derivation with the differentiation.

Hence, because of the linearity property for $(p,q)$-extensors it follows
that the $A$-directional derivative of $F$ at $X_{0}$ has the remarkable
property: for any $\alpha,\beta\in\mathbb{R}$ and $A,B\in\bigwedge^{p}V$
\begin{equation}
F_{\alpha A+\beta B}^{\prime}(X_{0})=\alpha F_{A}^{\prime}(X_{0})+\beta
F_{B}^{\prime}(X_{0}).   \label{6.3d}
\end{equation}

\begin{proposition}
Let $X:S\rightarrow\Lambda^{q}V$ be any $q$-vector function of a real
variable derivable at $\lambda_{0}\in S.$ Then, $X$ is differentiable at $%
\lambda_{0}$ and the differential of $X$ at $\lambda_{0}$ is $%
X_{\lambda_{0}}\in ext_{0}^{q}(V)$ given by
\begin{equation}
X_{\lambda_{0}}(\alpha)=X^{\prime}(\lambda_{0})\alpha,   \label{6.3e}
\end{equation}
where $X^{\prime}(\lambda_{0})$ is the derivative of $X$ at $\lambda_{0}.$
\end{proposition}

\begin{proof}
We only need to prove that
\begin{equation*}
\underset{\lambda\rightarrow\lambda_{0}}{\lim}\frac{X(\lambda)-X(\lambda
_{0})-X^{\prime}(\lambda_{0})(\lambda-\lambda_{0})}{\left| \lambda
-\lambda_{0}\right| }=0.
\end{equation*}

Since $X$ is derivable at $\lambda_{0},$ there is a $q$-vector function of
real variable, say $\xi_{\lambda_{0}},$ continuous at $\lambda_{0}$ such
that for all $\lambda\in S$
\begin{equation*}
X(\lambda)=X(\lambda_{0})+(\lambda-\lambda_{0})X^{\prime}(\lambda
_{0})+(\lambda-\lambda_{0})\xi_{\lambda_{0}}(\lambda),
\end{equation*}
where $\xi_{\lambda_{0}}(\lambda_{0})=0.$

Hence, it follows that for all $\lambda\neq\lambda_{0}$
\begin{equation*}
\frac{X(\lambda)-X(\lambda_{0})-X^{\prime}(\lambda_{0})(\lambda-\lambda_{0})%
}{\left| \lambda-\lambda_{0}\right| }=\pm\xi_{\lambda_{0}}(\lambda).
\end{equation*}
Thus, by taking limits for $\lambda\rightarrow\lambda_{0}$ on both sides, we
get the expected result.
\end{proof}

From eqs.(\ref{6.3c}) and (\ref{6.3e}), it should be noted that
the $\alpha $-directional derivative of $X$ at $\lambda_{0}$ is
given by
\begin{equation}
X_{\alpha}^{\prime}(\lambda_{0})=X^{\prime}(\lambda_{0})\alpha.
\label{6.3f}
\end{equation}

\subsubsection{Differentiation Rules}

Take two open subset of $\bigwedge^{p}V,$ say $\mathbf{\Omega}_{1}^{p}V$ and
$\mathbf{\Omega}_{2}^{p}V,$ such that $\mathbf{\Omega}_{1}^{p}V\cap \mathbf{%
\Omega}_{2}^{p}V\neq\emptyset.$

\begin{theorem}
Let $F:\mathbf{\Omega}_{1}^{p}V\rightarrow\bigwedge^{q}V$ and $G:\mathbf{%
\Omega}_{2}^{p}V\rightarrow\bigwedge^{q}V$ be two differentiable functions
at $X_{0}\in\mathbf{\Omega}_{1}^{p}V\cap\mathbf{\Omega}_{2}^{p}V$. Denote
the differentials of $F$ and $G$ at $X_{0}$ by $f_{X_{0}}$ and $g_{X_{0}},$
respectively.

The addition $F+G:\mathbf{\Omega}_{1}^{p}V\cap\mathbf{\Omega}%
_{2}^{p}V\rightarrow\Lambda^{q}V$ such that $(F+G)(X)=F(X)+G(X)$ and the
products $F\ast G:\mathbf{\Omega}_{1}^{p}V\cap\mathbf{\Omega}%
_{2}^{p}V\rightarrow \bigwedge V$ such that $(F\ast G)(X)=F(X)\ast G(X),$
where $\ast$ means either $(\wedge),$ $(\cdot),$ $(\lrcorner\llcorner)$ or $(
$\emph{Clifford product}$),$ are also differentiable functions at $X_{0}.$

The differential of $F+G$ at $X_{0}$ is $f_{X_{0}}+g_{X_{0}}$ and the
differentials of $F\ast G$ at $X_{0}$ are given by $A\mapsto
f_{X_{0}}(A)\ast G(X_{0})+F(X_{0})\ast g_{X_{0}}(A).$
\end{theorem}

\begin{proof}
We must prove that $s_{X_{0}}=f_{X_{0}}+g_{X_{0}}\in
ext_{p}^{q}(V)$ satisfies
\begin{equation*}
\underset{X\rightarrow X}{\lim}_{0}\frac{%
(F+G)(X)-(F+G)(X_{0})-s_{X_{0}}(X-X_{0})}{\left\| X-X_{0}\right\| }=0.
\end{equation*}
And, $p_{X_{0}}(A)=f_{X_{0}}(A)*G(X_{0})+F(X_{0})*g_{X_{0}}(A),$ more
general extensors over $V,$ verify

\begin{equation*}
\underset{X\rightarrow X}{\lim}_{0}\dfrac{%
(F*G)(X)-(F*G)(X_{0})-p_{X_{0}}(X-X_{0})}{\left\| X-X_{0}\right\| }=0.
\end{equation*}

Since $F$ and $G$ are differentiable at $X_{0},$ there are two $(p,q)$%
-functions $\varphi_{X_{0}}$ and $\psi_{X_{0}},$ continuous at $X_{0}$, such
that for all $X\in\mathbf{\Omega}_{1}^{p}V\cap\mathbf{\Omega}_{2}^{p}V$%
\begin{align*}
F(X) & =F(X_{0})+f_{X_{0}}(X-X_{0})+\left\| X-X_{0}\right\|
\varphi_{X_{0}}(X), \\
G(X) & =G(X_{0})+g_{X_{0}}(X-X_{0})+\left\| X-X_{0}\right\| \psi_{X_{0}}(X),
\end{align*}
where $\varphi_{X_{0}}(X_{0})=\psi_{X_{0}}(X_{0})=0$.

Hence, the following multivector identities which hold for all $X\neq X_{0}$
can be easily deduced
\begin{equation*}
\dfrac{(F+G)(X)-(F+G)(X_{0})-s_{X_{0}}(X-X_{0})}{\left\| X-X_{0}\right\| }%
=\varphi_{X_{0}}(X)+\psi_{X_{0}}(X)
\end{equation*}

and

\begin{align*}
& \dfrac{(F*G)(X)-(F*G)(X_{0})-p_{X_{0}}(X-X_{0})}{\left\| X-X_{0}\right\| }
\\
= & \varphi_{X_{0}}(X)*G(X_{0})+F(X_{0})*\psi_{X_{0}}(X) \\
& +\varphi_{X_{0}}(X)*g_{X_{0}}(X-X_{0})+f_{X_{0}}(X-X_{0})*\psi_{X_{0}}(X)
\\
& +\frac{f_{X_{0}}(X-X_{0})*g_{X_{0}}(X-X_{0})}{\left\| X-X_{0}\right\| }%
+\left\| X-X_{0}\right\| \varphi_{X_{0}}(X)*\psi_{X_{0}}(X).
\end{align*}

Now, by taking limits for $X\rightarrow X_{0}$ on both sides of these
multivector identities\footnote{%
For calculating some limits we have used an useful lemma. For any $f\in
ext_{p}^{q}(V)$ there exists a real number $M\geq0$ such that for every $%
X\in\bigwedge^{p}V:\left\| f(X)\right\| \leq M\left\| X\right\| .$}, we get
the required results.
\end{proof}

In accordance with eq.(\ref{6.3c}) all differentiation rule turns out to be
an $A$-directional derivation rule.

For the addition of two differentiable functions $F$ and $G$ we have
\begin{equation*}
(F+G)_{A}^{%
\prime}(X_{0})=(f_{X_{0}}+g_{X_{0}})(A)=f_{X_{0}}(A)+g_{X_{0}}(A),
\end{equation*}
i.e.,
\begin{equation}
(F+G)_{A}^{\prime}(X_{0})=F_{A}^{\prime}(X_{0})+G_{A}^{\prime}(X_{0}).
\label{6.4a}
\end{equation}

For the products $F*G$ we get
\begin{equation*}
(F*G)_{A}^{\prime}(X_{0})=f_{X_{0}}(A)*G(X_{0})+F(X_{0})*g_{X_{0}}(A),
\end{equation*}
i.e.,
\begin{equation}
(F*G)_{A}^{\prime}(X_{0})=F_{A}^{\prime}(X_{0})*G(X_{0})+F(X_{0})*G_{A}^{%
\prime}(X_{0}).   \label{6.4b}
\end{equation}

\begin{theorem}
Take an open subset of $\bigwedge^{p}V,$ say $\mathbf{\Omega}^{p}V$. Let $G:%
\mathbf{\Omega}^{p}V\rightarrow\bigwedge^{q}V$ and $F:\bigwedge
^{q}V\rightarrow\bigwedge^{r}V$ be two differentiable functions, the first
one at $X_{0}\in\mathbf{\Omega}^{p}V$ and the second one at $G(X_{0})\in
\bigwedge^{q}V.$ Denote by $g_{X_{0}}$ and $f_{G(X_{0})}$ the differentials
of $G$ at $X_{0}$ and of $F$ at $G(X_{0}),$ respectively.

The composition $F\circ G:\mathbf{\Omega}^{p}V\rightarrow\bigwedge^{r}V$
such that $F\circ G(X)=F(G(X))$ is also a differentiable function at $X_{0}.$
The differential of $F\circ G$ at $X_{0}$ is $f_{G(X_{0})}\circ g_{X_{0}}.$
\end{theorem}

\begin{proof}
We must prove that for $f_{G(X_{0})}\circ g_{X_{0}}\in ext_{p}^{r}(V)$ it
holds
\begin{equation*}
\underset{X\rightarrow X_{0}}{\lim}\frac{F\circ G(X)-F\circ
G(X_{0})-f_{G(X_{0})}\circ g_{X_{0}}(X-X_{0})}{\left\| X-X_{0}\right\| }=0.
\end{equation*}

Since $G$ is differentiable at $X_{0}$ and $F$ is differentiable at $%
G(X_{0}),$ there are a $(p,q)$-function $X\mapsto\psi_{X_{0}}(X)$ and a $%
(q,r)$-function $Y\mapsto\varphi_{G(X_{0})}(Y),$ the first one continuous at
$X_{0}$ and the second one continuous at $G(X_{0}),$ such that for all $X\in%
\mathbf{\Omega}^{p}V$ and $Y\in\bigwedge^{q}V$
\begin{align*}
G(X) & =G(X_{0})+g_{X_{0}}(X-X_{0})+\left\| X-X_{0}\right\| \psi_{X_{0}}(X),
\\
F(Y) & =F(G(X_{0}))+f_{G(X_{0})}(Y-G(X_{0}))+\left\| Y-G(X_{0})\right\|
\varphi_{G(X_{0})}(Y),
\end{align*}
where $\psi_{X_{0}}(X_{0})=0$ and $\varphi_{G(X_{0})}(G(X_{0}))=0.$

Hence, it follows easily a multivector identity which holds for all $X\neq
X_{0}$
\begin{align*}
& \frac{F\circ G(X)-F\circ G(X_{0})-f_{G(X_{0})}\circ g_{X_{0}}(X-X_{0})}{%
\left\| X-X_{0}\right\| } \\
& =f_{G(X_{0})}\circ\psi_{X_{0}}(X)+\frac{\left\| G(X)-G(X_{0})\right\| }{%
\left\| X-X_{0}\right\| }\varphi_{G(X_{0})}\circ G(X).
\end{align*}

Now, by taking limits for $X\rightarrow X_{0}$ on both sides, using the
equations: $\underset{X\rightarrow X_{0}}{\lim}f_{G(X_{0})}\circ%
\psi_{X_{0}}(X)=0$ and $\underset{X\rightarrow X_{0}}{\lim}\dfrac{\left\|
G(X)-G(X_{0})\right\| }{\left\| X-X_{0}\right\| }\varphi_{G(X_{0})}\circ
G(X)=0,$ we get the expected result.
\end{proof}

This \emph{chain rule} for differentiation turns out to be a chain rule for $%
A$-directional derivation.

For a differentiable $G$ at $X_{0}$ and a differentiable $F$ at $G(X_{0})$
we have
\begin{align}
(F\circ G)_{A}^{\prime}(X_{0}) & =f_{G(X_{0})}\circ
g_{X_{0}}(A)=f_{G(X_{0})}(G_{A}^{\prime}(X_{0})),  \notag \\
(F\circ G)_{A}^{\prime}(X_{0}) &
=F_{G_{A}^{\prime}(X_{0})}^{\prime}(G(X_{0})).   \label{6.4c}
\end{align}

We study now two very important particular cases of the general chain rule
for the $A$-directional derivation:

For $p>0,$ $q=0$ and $r>0,$ i.e., for the $A$-directional derivative of the
composition of $\Phi:\mathbf{\Omega}^{p}V\rightarrow\mathbb{R}$ with $X:%
\mathbb{R}\rightarrow\bigwedge^{r}V$ at $X_{0}\in\mathbf{\Omega}^{p}V,$ by
using eq.(\ref{6.4c}) and eq.(\ref{6.3f}), we have
\begin{align}
(X\circ\Phi)_{A}^{\prime}(X_{0}) & =X_{\Phi_{A}^{\prime}(X_{0})}^{\prime
}(\Phi(X_{0})),  \notag \\
(X\circ\Phi)_{A}^{\prime}(X_{0}) & =X^{\prime}(\Phi(X_{0}))\Phi_{A}^{\prime
}(X_{0}).   \label{6.4d}
\end{align}

For $p=0,$ $q>0$ and $r>0,$ i.e., for the derivative of the composition of $%
X:S\rightarrow\bigwedge^{q}V$ with $F:\bigwedge^{q}V\rightarrow\bigwedge
^{r}V$ at $\lambda_{0}\in S,$ by using eq.(\ref{6.3f}), eq.(\ref{6.4c}) and
eq.(\ref{6.3d}), we have
\begin{align}
(F\circ X)^{\prime}(\lambda_{0})\alpha & =(F\circ
X)_{\alpha}^{\prime}(\lambda_{0})=F_{X_{\alpha}^{\prime}(\lambda_{0})}^{%
\prime}(X(\lambda _{0})),  \notag \\
(F\circ X)^{\prime}(\lambda_{0})\alpha & =F_{X^{\prime}(\lambda_{0})\alpha
}^{\prime}(X(\lambda_{0})).   \label{6.4e}
\end{align}

\subsubsection{Derivatives}

Let $(\{e_{k}\},\{e^{k}\})$ be a pair of reciprocal bases of $V.$ Let $F:%
\mathbf{\Omega}^{p}V\rightarrow\bigwedge^{q}V$ be any differentiable
function at $X_{0}\in\mathbf{\Omega}^{p}V.$ Define the set $\mathbf{\Lambda }%
^{p}V=\{X\in\mathbf{\Omega}^{p}V$ $/$ $F$ is differentiable at $X\}\subseteq
\mathbf{\Omega}^{p}V.$

It follows that it must exist a well-defined function $F_{A}^{\prime }:%
\mathbf{\Lambda}^{p}V\rightarrow\bigwedge^{q}V$ such that $F_{A}^{\prime }(X)
$ equals the $A$-directional derivative of $F$ at each $X\in \mathbf{\Lambda}%
^{p}V.$ It is called the $A$-\emph{directional derivative function of} $F$.

Then, we can define exactly \emph{four} derivative-like functions for $F$,
namely, $*F^{\prime}:\mathbf{\Lambda}^{p}V\rightarrow\bigwedge^{q}V$ such
that
\begin{align}
\ast F^{\prime}(X) & =\frac{1}{p!}(e^{j_{1}}\wedge\ldots
e^{j_{p}})*F_{e_{j_{1}}\wedge\ldots e_{j_{p}}}^{\prime}(X)  \notag \\
& =\frac{1}{p!}(e_{j_{1}}\wedge\ldots e_{j_{p}})*F_{e^{j_{1}}\wedge\ldots
e^{j_{p}}}^{\prime}(X),   \label{6.5a}
\end{align}
where $*$ means either $(\wedge),$ $(\cdot),$ $(\lrcorner)$ or $($\emph{%
Clifford product}$).$

Whichever $*F^{\prime}$ is a well-defined function associated to $F,$ since $%
*F^{\prime}(X)$ are multivectors which do not depend on the choice of $%
(\{e_{k}\},\{e^{k}\}).$

We will call $\wedge F^{\prime},$ $\cdot F^{\prime},$ $\lrcorner F^{\prime}$
and $F^{\prime}$ (i.e., $*\equiv($\emph{Clifford product}$)$) respectively
the (generalized) \emph{curl,} \emph{scalar divergence,} \emph{left
contracted divergence }and \emph{gradient of }$F.$ Sometimes the gradient of
$F$ is called the \emph{standard derivative of }$F.$

On the \emph{real vector space of differentiable }$(p,q)$\emph{-functions
over }$V$ we can introduce exactly four derivative-like operators, namely, $%
F\mapsto\partial_{X}*F$ such that
\begin{equation}
\partial_{X}*F=*F^{\prime},   \label{6.5b}
\end{equation}
i.e., for every $X\in\mathbf{\Lambda}^{p}V$
\begin{align}
\partial_{X}*F(X) & =\frac{1}{p!}(e^{j_{1}}\wedge\ldots
e^{j_{p}})*F_{e_{j_{1}}\wedge\ldots e_{j_{p}}}^{\prime}(X)  \notag \\
& =\frac{1}{p!}(e_{j_{1}}\wedge\ldots e_{j_{p}})*F_{e^{j_{1}}\wedge\ldots
e^{j_{p}}}^{\prime}(X).   \label{6.5c}
\end{align}
The special cases $\partial_{X}\wedge$, $\partial_{X}\cdot,$ $\partial
_{X}\lrcorner$ and $\partial_{X}$ (i.e., $*\equiv($\emph{Clifford product}$)$%
) will be called respectively the (generalized) \emph{curl, scalar
divergence,} \emph{left contracted} \emph{divergence }and\emph{\ gradient
operator.}

For differentiable functions it is also possible to introduce a remarkable
operator denoted by $A\cdot\partial_{X},$ and defined as follows
\begin{align}
A\cdot\partial_{X}F(X) & =(A\cdot\frac{1}{p!}e^{j_{1}}\wedge\ldots
e^{j_{p}})F_{e_{j_{1}}\wedge\ldots e_{j_{p}}}^{\prime}(X),  \notag \\
& =(A\cdot\frac{1}{p!}e_{j_{1}}\wedge\ldots
e_{j_{p}})F_{e^{j_{1}}\wedge\ldots e^{j_{p}}}^{\prime}(X),   \label{6.5d}
\end{align}
i.e., by eq.(\ref{6.3d})
\begin{equation}
A\cdot\partial_{X}F(X)=F_{A}^{\prime}(X).   \label{6.5e}
\end{equation}
The operator $A\cdot\partial_{X}$ is called the $A$-directional derivative
operator. It maps $F$ $\rightarrow$ $F_{A}^{\prime},$ i.e., $A\cdot
\partial_{X}F=F_{A}^{\prime}.$

Now, we write out the property expressed by eq.(\ref{6.3d}) using the
operator $A\cdot\partial_{X}$. We have,
\begin{equation}
(\alpha A+\beta B)\cdot\partial_{X}F(X_{0})=\alpha
A\cdot\partial_{X}F(X_{0})+\beta B\cdot\partial_{X}F(X_{0}).   \label{6.5f}
\end{equation}
We have then a suggestive operator identity
\begin{equation}
(\alpha A+\beta B)\cdot\partial_{X}=\alpha A\cdot\partial_{X}+\beta
B\cdot\partial_{X}.   \label{6.5ff}
\end{equation}

We have also rules holding for the $A$-directional derivation of addition,
products and composition of differentiable functions, and eq.(\ref{6.4a}),
eq.(\ref{6.4b}) and eq.(\ref{6.4c}) can be written as:
\begin{equation}
A\cdot\partial_{X}(F+G)(X_{0})=A\cdot\partial_{X}F(X_{0})+A\cdot\partial
_{X}G(X_{0}),   \label{6.5g}
\end{equation}
i.e., $A\cdot\partial_{X}(F+G)=A\cdot\partial_{X}F+A\cdot\partial_{X}G$.
\begin{equation}
A\cdot\partial_{X}(F*G)(X_{0})=A\cdot%
\partial_{X}F(X_{0})*G(X_{0})+F(X_{0})*A\cdot\partial_{X}G(X_{0}),
\label{6.5h}
\end{equation}
i.e., $A\cdot\partial_{X}(F*G)=(A\cdot\partial_{X}F)*G+F*(A\cdot\partial
_{X}G).$

If $X\mapsto G(X)$ and $Y\mapsto F(Y),$ then
\begin{equation}
A\cdot\partial_{X}(F\circ G)(X_{0})=A\cdot\partial_{X}G(X_{0})\cdot
\partial_{Y}F(G(X_{0})).   \label{6.5i}
\end{equation}

\section{Conclusions}

We studied in detail the concepts of limit, continuity and differentiability
for multivector functions of a $p$-vector variable. Several types of
derivatives for these objects have been introduced as, e.g., the $A$%
-directional derivative and the generalized curl, divergence and gradient.
We saw that the concept of extensor plays a key role in the formulation of
the notion of differentiability, it implies the existence of the
differential extensor. We have proved the basic derivation rules for all
suitable products of multivector functions and for composition of
multivector functions. The generalization of these results towards a
formulation of a general theory of multivector functions of several
multivector variables can be done easily. The concept of multivector
derivatives has been first introduced in \cite{3}. We think that our
presentation is an improvement of that presentation, clearing many issues.

In the following paper about multivector \emph{functionals}, we will see
that the gradient-derivative plays a key role in the formulation of
derivation concepts for the so-called induced multivector
functionals.\bigskip

\textbf{Acknowledgement}: V. V. Fern\'{a}ndez is grateful to FAPESP for a
posdoctoral fellowship. W. A. Rodrigues Jr. is grateful to CNPq for a senior
research fellowship (contract 201560/82-8) and to the Department of
Mathematics of the University of Liverpool for the hospitality. Authors are
also grateful to Drs. P. Lounesto, I. Porteous, and J. Vaz, Jr. for their
interest on our research and useful discussions.

\end{document}